				\documentclass[11pt, dvipdfmx]{amsart}  
	
	\usepackage{graphicx, latexsym, amsmath, amsthm, amssymb, mystyle, color} 
	
\def\lk{\mathrm{lk}} 	\def\det{\mathrm{det}} 
\def\itr{\mathrm{int}} 
	
\def\sh{\:\sharp\:}
\def\ms{\bb-\bb}

\def\borderarray#1#2#3#4#5#6{%
\setbox0\hbox{$\begin{array}{#5}#6\end{array}$}
\setlength{\dimen1}{\wd0}\addtolength{\dimen1}{-#3}\addtolength{\dimen1}{-\arraycolsep}
\setlength{\dimen2}{\ht0}\addtolength{\dimen2}{-#4}
\setbox1\hbox{$\left#1\rule{\dimen1}{0pt}\rule{0pt}{\dimen2}\right#2$}
\setbox0\hbox{\raisebox{\dp0}{\box0}\kern-\dimen1\kern-5pt\raisebox{\dp1}{\box1}}
\vcenter{\box0}
}

	\title{\bf{	Alexander polynomials of simple-ribbon knots}}
	\author{Kengo Kishimoto, Tetsuo Shibuya, Tatsuya Tsukamoto, and Tsuneo Ishikawa}
	\thanks{This work was supported by JSPS KAKENHI Grant Number JP16K05162.}
	\date{}
\begin{document}
                         \baselineskip=16pt 
                         \maketitle  

	\bgnC \bf{Abstract} \endC

		This is a revised version of \cite{KSTI-OJM21}. We revised the diagrams of $10_{42}$,
		$10_{75}$ in Figure \ref{fig:exp_srkten}, and the values of $l$ of $F(t\,;m,p,l)$ for
		$10_{42}$, $10_{75}$ and the value of $p$ of $F(t\,;m,p,l)$ for $10_{99}$ in Table 1.
		In \cite{KST-JMSJ16}, we introduced special types of fusions, so called simple-ribbon 
		fusions on links. A knot obtained from the trivial knot by a finite sequence of 
		simple-ribbon fusions is called a simple-ribbon knot. Every ribbon knot with $\leq 9$ 
		crossings is a simple-ribbon knot. In this paper, we give a formula for the Alexander 
		polynomials of simple-ribbon knots. Using the formula, we determine if a knot with 
		$10$ crossings is a simple-ribbon knot. Every simple-ribbon fusion can be realized by 
		``elementary" simple-ribbon fusions. We call a knot an $m$-simple-ribbon knot if the 
		knot is obtained from the trivial knot by a finite sequence of elementary 
		$m$-simple-ribbon fusions for a fixed positive integer $m$. We provide a condition 
		for a simple-ribbon knot to be both of an $m$-simple-ribbon knot and an 
		$n$-simple-ribbon knot for 	positive integers $m$ and $n$. 
%
%-----------------------------------------------------------------------------------------------------------------------------------
%
 
						\section{Introduction}\label{sec:intro}
									
		Knots and links are assumed to be ordered and oriented, and they are considered 
		up to ambient isotopy in an oriented $3$-sphere $S^3$. In \cite{KST-JMSJ16}, we 
		introduced special types of fusions, so called simple-ribbon fusions.
		A ($m$-){\it ribbon fusion} on a link $\ell$ is an $m$-{\it fusion} (\cite[Definition 13.1.1]
		{AK-book96}) on the split union of $\ell$ and an $m$-component trivial link $\mO$
		such that each component of $\mO$ is attached to a component of $\ell$ by a single 
		band. Note that any knot obtained from the trivial knot by a finite sequence of ribbon 
		fusions is a {\it ribbon knot} (\cite[Definition 13.1.9]{AK-book96}), and that any ribbon 
		knot can be obtained from the trivial knot by ribbon fusions. Here we only define an 
		elementary simple-ribbon fusion. A general simple-ribbon fusion can be realized by 
		elementary simple-ribbon fusions. Refer \cite{KST-JMSJ16} for precise definition.
\pvb
		Let $\ell$ be a link and $\mO=O_1\ccc O_m$ the $m$-component trivial link which is 
		split from $\ell$. Let $\mD=D_1\ccc D_m$ be a disjoint union of non-singular disks with 
		$\ptl D_i=O_i$ and $D_i\cap\ell=\emptyset$ ($i=1,\cdots, m$), and let $\mB=B_1\ccc B_m$ 
		be a disjoint union of disks for an $m$-fusion, called {\it bands}, on the split union of 
		$\ell$ and $\mO$ satisfying the following (see Figure \ref{fig:itr_explsrk} for example):\pvb
\bgnI
		\item[(i)] 	$B_i\cap\ell=\ptl B_i\cap \ell =\{$ a single arc $\}$;\pvb
		\item[(ii)] 	$B_i\cap\mO=\ptl B_i\cap O_i =\{$ a single arc $\}$; and\pvb
		\item[(iii)] 	$B_i\cap{\rm int}\: \mD=B_i\cap{\rm int}\: D_{i+1}=\{$ a single arc of ribbon 
					type $\}$.	
\endI
\bgnF  		
		{\iclg[scale=1, bb=0 0 422 157]{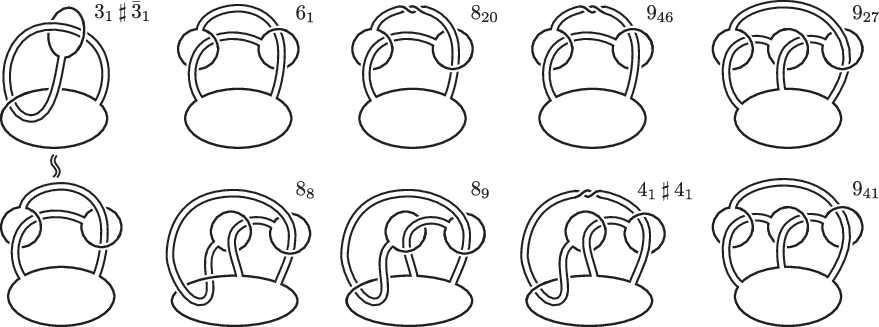}} 
		\captn{ribbon knots with less than or equal to nine crossings}\label{fig:itr_explsrk}	\endF
\pvc	
		Let $L$ be a link obtained from the split union of $\ell$ and $\mO$ by the $m$-fusion along 
		$\mB$, i.e., $L=(\ell\cup\mO\cup\ptl\mB)-{\rm int}(\mB\cap\ell)-{\rm int}(\mB\cap\mO)$. Then 
		we say that $L$ is obtained from $\ell$ by an {\it elementary} ($m$-){\it simple-ribbon fusion}  
		or an {\it elementary} ($m$-)SR-fusion ({\it with respect to} $\mD\cup\mB$). If a knot $K$ 
		is obtained from the trivial knot $O$ by a finite sequence of elementary SR-fusions, then we 
		call $K$ a {\it simple-ribbon knot} (or an SR-knot).  Since an elementary SR-fusion is a ribbon 
		fusion, any SR-knot is a ribbon knot. We also call the trivial knot an SR-knot. As illustrated 
		in Figure \ref{fig:itr_explsrk}, all the ribbon knots with $\leq 9$ crossings are SR-knots. 
\npg
		Let $\dot{D}_i$ and $\dot{B}_i$ be disks and $f: \cup_i \lpa \dot{D}_i \cup \dot{B}_i
		\rpa\to S^3$ an immersion such that $f(\dot{D}_i)=D_i$ and $f(\dot{B}_i)=B_i$. 
		We denote the arc of ${\rm int}\:D_i \cap B_{i-1}$ by $\a_i$ and let $B_{i,1}$ and $B_{i,2}$ 
		be the subdisks of $B_i$ such that $B_{i,1}\cup B_{i,2}=B_i$, $B_{i,1}\cap B_{i,2}=\a_{i+1}$, 
		and $B_{i,1}\cap \ptl D_i\neq\emptyset$. Take a point $b_i$ on ${\rm int}\:\a_i$ ($i=1$, 
		$\ldots$, $m$) and an arc $\b_i$ on $D_i\cup B_{i,1}$ so that $\b_i\cap(\a_i\cup\a_{i+1})
		=\ptl\b_i=b_i\cup b_{i+1}$ and oriented from $b_{i+1}$ to $b_i$ 
		(see Figure \ref{fig:itr_attknot}). 
		Then $\b=\cup_i \b_i$ is an oriented simple loop and we call $\b$ an {\it attendant knot} 
		of $\mD\cup\mB$. Moreover, we denote the pre-images of $\a_i$ (resp. $b_i$) on 
		$\dot{D}_i$ and $\dot{B}_{i-1}$ by $\dot{\a}_i$ and $\ddot{\a}_i$ (resp.  
		$\dot{b}_i$ and $\ddot{b}_i$), respectively.  
\bgnF  		
		{\iclg[scale=1, bb=0 0 423 124]{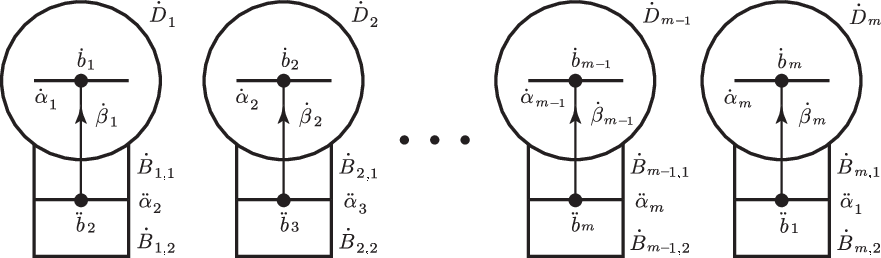}} \captn{}\label{fig:itr_attknot}	\endF
\pva
		$\mD\cup\mB$ is oriented so that induced orientations on 
		boundaries are compatible with the orientation of $\ell$. Then we can see that each 
		band $B_i$ intersects with $D_{i+1}$ in two ways, i.e. when we pass through $\a_{i+1}$ 
		from $B_{i,2}$ to $B_{i,1}$, 	we pass through $D_{i+1}$ either from the negative side to 
		the positive side of $D_{i+1}$, or from the positive side to the negative side of $D_{i+1}$. 
		In the former and latter cases, we say that $B_i$ is {\it positive} and {\it negative}, 
		respectively. Then we have the following.	
		
\bgnT	\label{thm:alexpoly}
		Let $K$ be a knot obtained from a knot $k$ by an elementary $m$-SR-fusion  
		with an attendant knot $\b$ and with $p$ positive bands.
		Let $l=\lk(\b, k)$ and $\v(t\,;m,p,l)=(1-t)^m - t^{\,l} (-t)^p$. Then we have the following.
		\begin{equation}\label{eq1a}
		\De_K(t)\doteq\De_k(t)\,\v(t\,;m,p,l)\,\v(t^{-1}\,;m,p,l)\end{equation}
		\endT
\nid		
		{\bf Remark.}	 
		We can also write $\De_K(t)$ as $\De_k(t)\,\v(t\,;m,p,l)\,\v(t\,;m,m\!-\!p,-l)$, i.e.
		\begin{equation}\label{eq1b}
		\De_K(t)\doteq\De_k(t)\,\{(1-t)^m - t^{\,l} (-t)^p\}\,\{(1-t)^m - t^{-l} (-t)^{m-p}\}\end{equation}
		
\bgnO	\label{cor:alexpoly}
		Let $K$ be a knot obtained from a knot $k$ by a finite sequence of elementary 
		SR-fusions, i.e., there exists a finite sequence $k=K_0$, $K_1$, $\ldots$, $K_N=K$ 
		of knots such that $K_i$ is obtained from $K_{i-1}$ by an elementary $m_i$-SR-fusion 
		with an attendant knot $\b_i$ and with $p_i$ positive bands $(i=1,\ldots, N)$. 
		Let $l_i=\lk(\b_i, K_{i-1})$ and $\v(t\,;m_i,p_i,l_i)=(1-t)^{m_i} - t^{\,l_i} (-t)^{p_i}$. 	
		Then we have the following.
		\begin{equation}\label{eq1c}
		\De_K(t)\doteq\De_k(t)\prod_{i=1}^N \v(t\,;m_i,p_i,l_i)\, \v(t^{-1}\,;m_i,p_i,l_i)\end{equation}
		\endO
	
		As mentioned in the beginning, all the ribbon knots with $\leq 9$ crossings are SR-knots. 
		Using Corollary \ref{cor:alexpoly}, we can determine if a ribbon knot with $10$ crossings 
		is an SR-knot. To do this, we use a value derived from the Alexander polynomial.
		For a knot $K$, let $\De'_K(t)$ be the polynomial such that $\De'_K(t)\doteq\De_K(t)$ 
		and $\De'_K(0)\neq 0$. Then define $\d_2(K)$ as $0$ if $|\De'_K(2)|=0$ and as the
		largest odd factor of $|\De'_K(2)|$ if $|\De'_K(2)|\neq 0$.
		Note that if $K$ is a simple-ribbon knot, then $\d_2(K)$ is a product of the 
		integers of the form $2^s\pm 1$ $(s=0, 1, 2, \ldots)$ from Corollary \ref{cor:alexpoly}. 

\bgnL	\label{lem:d2}
		If $K$ is a simple-ribbon knot such that $\d_2(K)=1$, then we have the following
		for a non-negative integer $n$. \begin{equation}\label{eq1d} 
		\De'_K(t)=1~~\mbox{\rm{or}}~~(1-6t+11t^2-6t^3+t^4)^n \end{equation}
		\endLR
		Since $K$ is a simple-ribbon knot, we have the following from Corollary \ref{cor:alexpoly}, 
		where $N$ $(\geq 1)$, $m_i$ $(\geq 1)$, $p_i$ $(0\leq p_i\leq m_i)$, and $l_i$ are 
		integers $(i=1, 2, \ldots, N)$. \pvbn
		$\De_K(t)\doteq\prod_{i=1}^N
		\{(1-t)^{m_i} - t^{\,l_i} (-t)^{p_i}\}\,\{(1-t)^{m_i} - t^{-l_i} (-t)^{m_i-p_i}\}$\pbg
		$\doteq \prod_{i=1}^N
		\{t^{p_i+l_i}+(-1)^{m_i-(p_i+1)}(t-1)^{m_i}\}\,\{t^{m_i-(p_i+l_i)}+(-1)^{p_i+1}(t-1)^{m_i}\}$
\pvc
		Let $g_i(t)=t^{p_i+l_i}+(-1)^{m_i-(p_i+1)}(t-1)^{m_i}$ and
		$h_i(t)=t^{m_i-(p_i+l_i)}+(-1)^{p_i+1}(t-1)^{m_i}$. Then we have that
		$\De'_K(2)=2^q\prod_{i=1}^Ng_i(2)h_i(2)$ for an integer $q$.
		Since $\d_2(K)=1$, each of	$|g_i(2)|$ and $|h_i(2)|$ is a power of $2$, and thus 
		$2^{-1}=|2^{-1}-1|$, $2=2^0+1$, or $1=2^1-1$	($i=1, 2, \ldots, N$). Thus, each of 
		$p_i+l_i$ and $m_i-(p_i+l_i)$ is $-1$, $0$, or $1$ for each $i$, and hence 
		$m_i=(p_i+l_i)+(m_i-(p_i+l_i))$ is $1$ or $2$, since $m_i>0$. Therefore we have that 
		$(g_i(2), h_i(2), m_i)=(2^0+1,2^1-1,1)$, $(2^1-1,2^0+1,1)$, or $(2^1-1, 2^1-1, 2)$. 
		In the first two cases and the last case, we have that
		$g_i(t)h_i(t)=\{t^0+(t-1)\}\{t^1-(t-1)\}=t$ and 
		$g_i(t)h_i(t)=\{t-(t-1)^2\}^2=1-6t+11t^2-6t^3+t^4$, respectively.
		Hence we obtain the conclusion.
		\endR
		
\bgnP	\label{prp:tenx}
		Among the $16$ ribbon knots with $10$ crossings, $10_{42}$, $10_{75}$, 
		$10_{87}$, $10_{99}$, $10_{129}$, $10_{137}$, $10_{140}$, $10_{153}$, and 
		$10_{155}$ are simple-ribbon knots and $10_{3}$, $10_{22}$, $10_{35}$, $10_{48}$, 
		$10_{123}$, $5_1\sh 5_1^*$, and $5_2\sh 5_2^*$ are not.
		\endPR
		The former statement is from Figure \ref{fig:exp_srkten}. To show the latter statement, 
		we consider $\d_2$ for each knot. Since $\d_2(10_{22})=11$, $\d_2(10_{48})=7\prd 13
		=1\prd 91$, and $\d_2(5_1\sh 5_1^*)=11\prd 11=1\prd 121$ from Table \ref{tab:alex} 
		and none of $11$, $13$, $91$, and $121$ is $2^s\pm 1$ for a non-negative integer 
		$s$, we know that these $3$ knots are not simple-ribbon knots. For the other $4$ knots,
		we have that $\d_2(10_{3})=\d_2(10_{35})=\d_2(10_{123})=\d_2(5_2\sh 5_2^*)=1$, and 
		the following from Table \ref{tab:alex}. Hence we know that they are not 
		simple-ribbon knots from Lemma \ref{lem:d2}.
\pcc		$\De'_{10_{3}}(t)=   6- 13t+ 6t^2$, 
\hsC	$\De'_{10_{35}}(t)= 2- 12t+ 21t^2- 12t^3+2t^4,$
\pcc		$\De'_{10_{123}}(t)= (1- 3t+3t^2- 3t^3+t^4)^2, $
\hsg		$\De'_{5_2\sh 5_2^*}(t)= 4 - 12 t+17 t^2- 12 t^3+4t^4$
		\endR
\bgnF  		
		{\iclg[scale=1, bb=0 0 468 177]{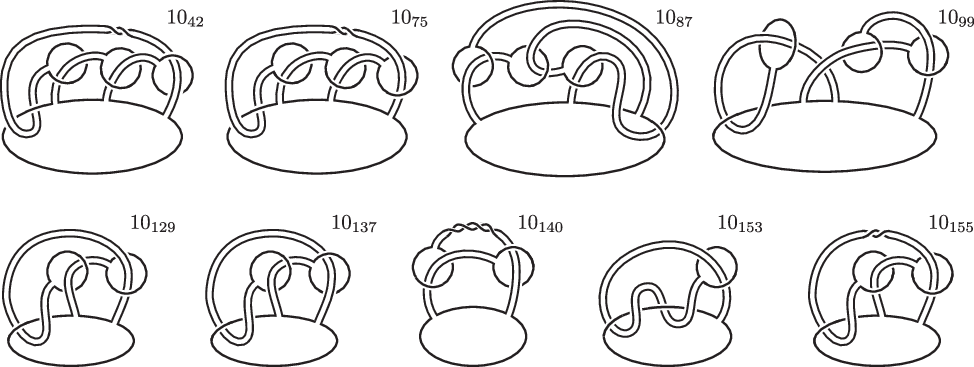}} \captn{}\label{fig:exp_srkten}	\endF
\pva
		Note that the above proof of Proposition \ref{prp:tenx} implies that for any ribbon knot
		$K$ with $\leq 10$ crossings, if $\De_K(t)$ can be written as equation (\ref{eq1c}), then
		$K$ is a simple-ribbon knot. However, it does not hold in general.		
		
\bgnT	\label{thm:nsrk}
		For any polynomial $\De(t)=\prod_{i=1}^N \v(t\,;m_i,p_i,l_i)\, \v(t^{-1}\,;m_i,p_i,l_i)$, 
		there exists a ribbon knot whose Alexander polynomial is $\De(t)$ and which is 
		not a simple-ribbon knot. 
		\endT

		If an SR-knot is obtained from the trivial knot by a finite sequence of elementary 
		$m$-SR-fusions for a fixed positive integer $m$, then we call the SR-knot $m$-SR-knot. 
		For example, $8_9$ is a $2$-SR-knot and $3_1\sh 3_1^*$ is a $1$-SR-knot and also 
		a $2$-SR-knot as we can see in Figure \ref{fig:itr_explsrk}. It is natural to ask if there 
		exists a simple-ribbon knot which is an $m$-SR-knot and also an $n$-SR-knot for 
		distinct positive integers $m$ and $n$ other than $3_1\sh 3_1^*$. We give a partial 
		answer to this question using equation (\ref{eq1c}). Let $m$ be a positive integer 
		and $\mK_m$ the set of non-trivial $m$-SR-knots. Then we have the following. 
		
\bgnT	\label{thm:msrk}
		If $\mK_m\cap\mK_n\neq\emptyset$ for positive integers $m$ and $n$ with $m>n$,
		then we have either that $(m, n)=(3,1)$, $(3,2)$, or $(2n,n)$.
		\endT

%-------------------------------------------------------------------------------------------------------------------------------------------------
%
%
%
%
%
%
%
%
%
%
%
%
%
%
%
%
%
%
%
%
%
%
%
%
%
%
%
%
%
%
%-------------------------------------------------------------------------------------------------------------------------------------------------

		\section{Proofs of Theorem \ref{thm:alexpoly} and Theorem \ref{thm:nsrk}}\label{sec:thmap}

		Let $K$ be a knot obtained from a knot $k$ by an elementary $m$-SR-fusion with 
		respect to $\mD\cup\mB$ with its attendant knot $\b$. Let $F$ be a Seifert surface 
		for $k$. Here we may take $F$ so that $F\cap\mD=\emptyset$. 
		Let $\mC=F\cup(\mD\cup\mB)$. We first transform $\mC$ into ``standard" position 
		and construct a Seifert surface $F_K$ for $K$ from $\mC$ in standard position. 
		Then, we calculate $\De_K(t)$ using $F_K$. 
\pvc
		We may take $F$ so that the intersections with $\mD\cup\mB$ are only arcs of 
		$\itr F$ and $\mB$. Then we divide the set of singularities of $\itr F\cap B_i$ into two: 
		one which consists of $\itr F\cap B_{i,1}$, denoted by $\mS_i$, and the other which 
		consists of $\itr F\cap B_{i,2}$, denoted by $\mT_i$. Thus the set of singularities of 
		$\mC$ is $\cup_i\a_i\cup\cup_i (\mS_i\cup\mT_i)$. We say that $\mC$ is {\it 
		in standard position} if 
		$\mS_1\cup\cdots\cup\mS_{m-1}=\emptyset$ and 
		$\mT_1\cup\cdots\cup\mT_m=\emptyset$ (see Figure \ref{fig:stp_expl} for example). 
		To transform $\mC$ into standard position, we need the following three transformations. 
		Here note that each transformation changes neither $m$, $p$, nor the knot type of $\b$. 
\pvcn
		{\bf Sliding a disk along a band} : Deforming $D_{i+1}$ by deformation retraction into
		a regular neighborhood of $B_i$ and slide $D_{i+1}$  along $B_i$ toward $D_i$. 
		Here $B_{i+1}$ follows $D_{i+1}$ (see Figure \ref{fig:stp_slide} for example). 
		We allow $D_{i+1}\cup B_{i+1}$ to pass through $F$. Then an additional intersection 
		of $B_{i+1}$ and $F$ is created.	
\pvcn
		{\bf Winding a band along $k$} : Winding $B_{i}$ along $k=\ptl F$ in a regular 
		neighborhood of $B_{i}\cap k$ either from negative side to positive side or from 
		positive side to negative side (see Figure \ref{fig:stp_twist} for example). Here an 
		additional intersection of $B_{i}$ and $F$ is created.
\pvcn	
		{\bf Tubing $F$} : Removing two disks $\d_1$ and $\d_2$ from $\itr\, F$ and attatch
		an annulus $S^1\prd[1,2]$ so that $S^1\prd \{i\}=\ptl\d_{i}$ $(i=1,2)$ and the result
		is orientable (see Figure \ref{fig:stp_tubing} for example).

\bgnP	\label{prop:stp}
		Let $K$ be a knot obtained from a knot $k$ by an elementary $m$-SR-fusion with 
		respect to $\mD\cup\mB$ with its attendant knot $\b$. Let $F$ be a Seifert surface 
		for $k$ such that $F\cap\mD=\emptyset$ and let $\mC\!=\!F\cup(\mD\cup\mB)$. 
		Then we may transform $\mC$ into standard position by sliding a disk along a band, 
		winding a band along $k$, and tubing $F$.
		\endPR
		First if $\mS_1\cup\cdots\cup\mS_{m-1}\neq\emptyset$, then take the smallest index 
		$i$ such that $\mS_i\neq\emptyset$ and slide $D_{i+1}$ along $B_i$ just next to $D_i$ 
		so that $\mS_i=\emptyset$ (See Figure \ref{fig:stp_slide} for example). Then slide 
		$D_{j+1}$ along $B_j$ inductively just next to $D_j$ so that $\mS_j=\emptyset$ 
		$(j=i+1,\ldots, m-1)$.
\vsmb
\bgnF  		
		{\iclg[scale=1, bb=0 0 426 158]{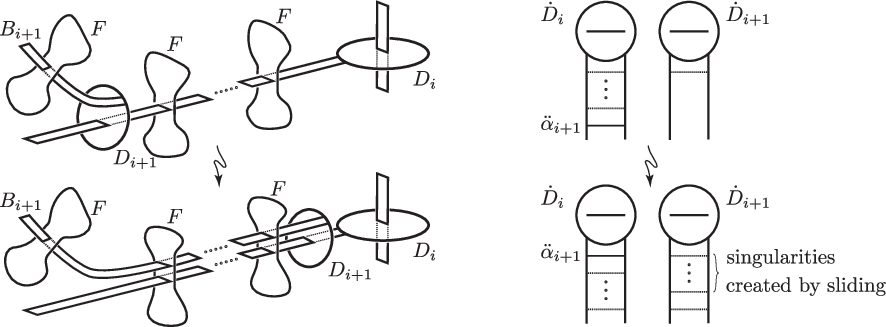}}\vsmc\captn{}\label{fig:stp_slide}	\endF
\par
\vsmb	Next if $\mT_1\cup\cdots\cup\mT_m\neq\emptyset$, then take an arbitrary $\mT_i\neq
		\emptyset$ and let $t_1$, $\ldots$, $t_p$ be its singularities which are placed close to 
		$B_i\cap k$ on $B_i$ in this order. Assume that $B_i$ is oriented as from $B_i\cap k$
		towards $B_i\cap D_i$ and let $\s(t_j)$ be the signed intersection number of $B_i$ 
		and $F$ at $t_j$. First wind $B_i$ along $k$ depending on $\s(t_j)$ $(j=1,\ldots, p)$. 
		If $\s(t_j)=1$ 	(resp. $-1$), then wind $B_i$ along $k=\ptl F$ from negative side to positive 
		side (resp. from positive side to negative side) as illustrated in Figure \ref{fig:stp_twist}. 
		Here we make these transformations from $j=1$ to $j=p$ in this order, and notice that 
		each transformation creates a new intersection $t'_j$ with $\s(t'_j)=-\s(t_j)$. Then make 
		a tubing $F$ so to erase $t_j$ and $t'_j$ from $j=1$ to $j=p$ in this order as illustrated 
		in Figure \ref{fig:stp_tubing}, and now $\mC$ is in standard position. \endR
\vsmc
\bgnF  		
		{\iclg[scale=.95, bb=0 0 419 133]{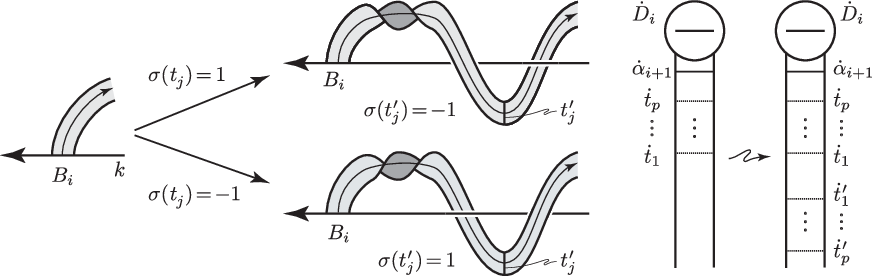}}\vsmb\captn{}\label{fig:stp_twist}	\endF
\vsme
\bgnF  		
		{\iclg[scale=1, bb=0 0 380 93]{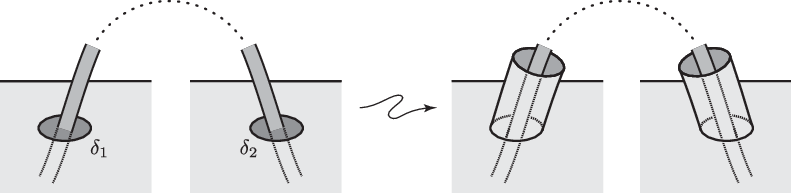}}\vsmb\captn{}\label{fig:stp_tubing}	\endF
%
%------------------------------------------------------------------------------------------------
%
\npg
\nid		{\bf\it Proof of Theorem \ref{thm:alexpoly}.}
		Let $F$ be a Seifert surface for $k$ such that $F\cap\mD=\emptyset$ and let 
		$\mC=F\cup(\mD\cup\mB)$. Here we may assume that $\mC$ is in standard 
		position from Proposition \ref{prop:stp}. Thus the set of singularities of $\mC$ is 
		$\cup_i\a_i\cup\mS_m$. Erase $\cup_i\a_i$ and $\mS_m$ to have a Seifert surface 
		$F_K$ for $K$ by orientation preserving cut and deformation as illustrated in the 
		second left of Figure \ref{fig:alx_basis_bd} and Figure \ref{fig:alx_basis_ls}, 
		respectively (see Figure \ref{fig:exp_basis} for example of $F_K$).
\bgnF  		
		{\iclg[scale=1, bb=0 0 433 165]{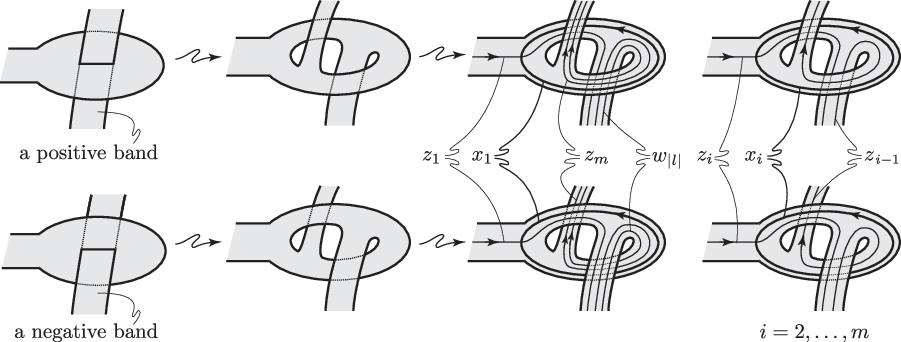}} \captn{}\label{fig:alx_basis_bd}\endF
\bgnF  		
		{\iclg[scale=1, bb=0 0 431 167]{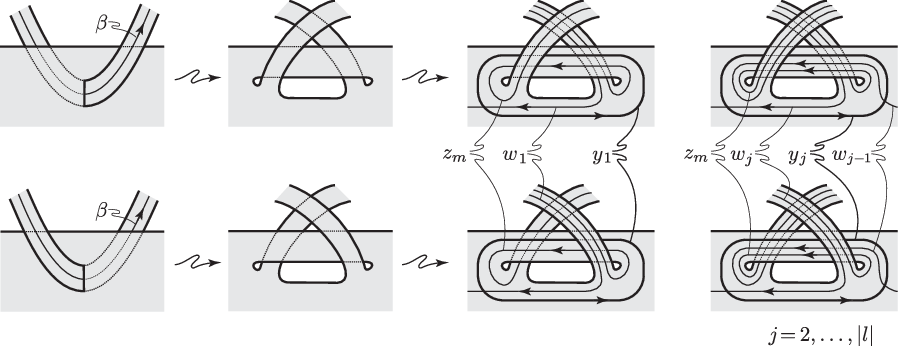}} \captn{}\label{fig:alx_basis_ls}	\endF
\par
		Take a basis $x_1$, $\ldots$, $x_m$, $y_1$, $\ldots$, $y_{|l|}$, $z_1$, $\ldots$, $z_m$, 
		$w_1$, $\ldots$, $w_{|l|}$, $u_1$, $\cdots$, $u_g$ of $H_1(F_K; \bZ)$ as illustrated in 
		Figure \ref{fig:alx_basis_bd} and Figure \ref{fig:alx_basis_ls} (see Figure \ref{fig:exp_basis} 
		for example), where $u_1$, $\cdots$, $u_g$ is a basis of $H_1(F; \bZ)$. 
		Then we have the following Seifert matrix $M$ with respect to the basis.
	\[M\ba=\ba
				\left( {\renewcommand\arraystretch{2} \begin{array}{c|c|c} 
				\bb O_{(m+|l|)\prd(m+|l|)} & P_{(m+|l|)\prd(m+|l|)} & O_{(m+|l|)\prd\,g\,}\bb 
		\\ \hline
				\bb Q_{(m+|l|)\prd(m+|l|)} & * & * \bb 
		\\ \hline 
				\bb O_{\,g\,\prd(m+|l|)} & * & M'  \bb \end{array}}\right)
		\ba=\ba
				\left( {\renewcommand\arraystretch{2} \begin{array}{c|c|c} 
				\bc O_{(m+|l|)\prd(m+|l|)} & 
			\begin{array}{c|c} \bc P^1_{m\prd m} 	& P^2_{m\prd |l|}\bc \\ \hline 
							\bc P^3_{|l|\prd m} 	& P^4_{|l|\prd |l|}\bc \end{array}
												& O_{(m+|l|)\prd\,g\,}\bc
		\\ \hline
			\begin{array}{c|c}	\bc Q^1_{m\prd m} & Q^2_{m\prd |l|} \bc\\ \hline 
							\bc Q^3_{|l|\prd m} & Q^4_{|l|\prd |l|} \bc\end{array} & * & * 
		\\ \hline
							\bc O_{\,g\,\prd(m+|l|)} 	  & * & M'  \end{array}}\right),\]

%
%------------------------------------------------------------------------------------------------
%
\npg
\nid
		where $M'$ is a Seifert matrix for $k$, $O_{s\times t}$ is the $s\times t$ zero matrix,
\pvan
		$P^1_{m\prd m}=(p^1_{ij})$ is an $m\prd m$ matrix with
		$p^1_{ij}=\mathrm{lk}(x_i,z_j^+)=\left\{\begin{array}{ll}
	 		\bc-(\e_i+1)/2 	& 	\text{if $i=j$} \\
			~\:~\:~\:~\:~\e_i 	& 	\text{if $i=1$ and $j=m$}, \\ &
								\text{or $2\leq i\leq m$ and $j=i-1$}\\
			~\:~\:~\:~\:~ 0 	& 	\text{otherwise,}\end{array}\right. $
\\
		$P^2_{m\prd |l|}=(p^2_{ij})$ is an $m\prd |l|$ matrix with
		$p^2_{ij}=\mathrm{lk}(x_i,w_j^+)=\left\{\begin{array}{ll}
	 	~\:~\:~\:~\:~ \e_1 ~\:~\:~\:~\:~\:~ 	& 	\text{if $i=1$ and $j=|l|$} \\
		~\:~\:~\:~\:~ 0  ~\:~\:~\:~\:~\:~ 	&	\text{otherwise,}\end{array}\right. $
\\
		$P^3_{|l|\prd m}=(p^3_{ij})$ is an $|l|\prd m$ matrix with
		$p^3_{ij}=\mathrm{lk}(y_i,z_j^+)~=\left\{\begin{array}{ll}
	 	~\:~\:~\:~\:~ \e ~\:~\:~\:~\:~\:~\:~ 	& 	\text{if $j=m$} \\
		~\:~\:~\:~\:~ 0 ~\:~\:~\:~\:~\:~\:~ 	& 	\text{otherwise,}\end{array}\right. $		
\\
		$P^4_{|l|\prd |l|}=(p^4_{ij})$ is an $|l|\prd |l|$ matrix with
		$p^4_{ij}=\mathrm{lk}(y_i,w_j^+)=\left\{\begin{array}{ll}
	 	(\e+1)/2\:~\:~ 	& \text{if $i=j$} \\
		(\e-1)/2\:~\:~		& \text{if $2\leq i\leq |l|$ and $j=i-1$}\\
		~\:~\:~\:~\:~ 0 	& \text{otherwise,}\end{array}\right. $	
\pvcn
		if $l\neq 0$, and $P_{m\prd m}=P^1_{m\prd m}$ if $l=0$,	
\\
		$Q^1_{m\prd m}=(q^1_{ij})$ is an $m\prd m$ matrix with
		$q^1_{ij}=\mathrm{lk}(z_i,x_j^+)=\left\{\begin{array}{ll}
	 	\bc-(\e_i-1)/2 	& 	\text{if $i=j$} \\
		~\:~\:~\:~\:~\e_i 	& 	\text{if $i=m$ and $j=1$}, \\ & 
							\text{or $1\leq i\leq m-1$ and $j=i+1$}\\
		~\:~\:~\:~\:~ 0 	& 	\text{otherwise,}\end{array}\right. $
\\
		$Q^2_{m\prd |l|}=(q^2_{ij})$ is an $m\prd |l|$ matrix with
		$q^2_{ij}=\mathrm{lk}(z_i,y_j^+)~=\left\{\begin{array}{ll}
	 	~\:~\:~\:~\:~ \e ~\:~\:~\:~\:~\:~\:~ 	& 	\text{if $i=m$} \\
		~\:~\:~\:~\:~ 0 ~\:~\:~\:~\:~\:~\:~ 	& 	\text{otherwise,}\end{array}\right. $
\\
		$Q^3_{|l|\prd m}=(q^3_{ij})$ is an $|l|\prd m$ matrix with
		$q^3_{ij}=\mathrm{lk}(w_i,x_j^+)=\left\{\begin{array}{ll}
	 	~\:~\:~\:~\:~ \e_1 ~\:~\:~\:~\:~\:~ 	& 	\text{if $i=|l|$ and $j=1$} \\
		~\:~\:~\:~\:~ 0  ~\:~\:~\:~\:~\:~ 	&	\text{otherwise, and}\end{array}\right. $		
\\
		$Q^4_{|l|\prd |l|}=(q^4_{ij})$ is an $|l|\prd |l|$ matrix with
		$q^4_{ij}=\mathrm{lk}(w_i,y_j^+)=\left\{\begin{array}{ll}
	 	(\e-1)/2\:~\:~ 		& 	\text{if $i=j$} \\
		(\e+1)/2\:~\:~		& 	\text{if $1\leq i\leq |l|-1$ and $j=i+1$}\\
		~\:~\:~\:~\:~\:~ 0 	& 	\text{otherwise,}\end{array}\right. $		
\pvcn
		if $l\neq 0$, and $Q_{m\prd m}=Q^1_{m\prd m}$ if $l=0$,	and
		$\e=\left\{\begin{array}{rl} 1 	& \text{if $l$ is positive}\\ 
						   	    -1  	& \text{if $l$ is negative}\end{array}\right. $,
		$\e_i=\left\{\begin{array}{rl} 1 	& \text{if $B_i$ is positive}\\ 
						   	       -1  	& \text{if $B_i$ is negative}\end{array}\right. $
\pvan
		$(i=1, \ldots, m)$.	Letting $a=\Frc{\e+1}{2}$, $b=\Frc{\e-1}{2}$, 
		$a_i=\Frc{\e_i+1}{2}$, and $b_i=\Frc{\e_i-1}{2}$, we have the following.
%
%------------------------------------------------------------------------------------------------
%
\pvcn\hspace{-3mm}{\small
		$P\ba=\ba\bordermatrix{
	   	&\ba z^+_1\ba&\ba z^+_2\ba&\ba\cdots\ba&\ba z^+_{m-1}\ba&\ba z_m\ba 
		&\ba w^+_1\ba&\ba w^+_2\ba&\ba\cdots\ba&\ba w^+_{|l|-1}\bb&\bb w^+_{|l|}\ba\cr
	x_1 	&-a_1\bd&&&&\bd\e_1\bd&&&&&\bd\e_1\bd\cr
	x_2 	&\e_2\bd&\bd-a_2\bd&&&&&&&&\cr
 ~\vdots &&\bd\ddots\bd&\bd\ddots\bd&&&&&&&\cr 
  x_{m-1} &&&\bd\ddots\bd&\bd-a_{m-1}\bd&&&&&&\cr
	x_m &&&&\bd\e_m\bd&\bd-a_m\bd&&&&&\cr
	y_1 	&&&&& \e &\bc a\bc&&&&\cr
	y_2 	&&&&& \e &\bc b\bc&\bc a\bc&&&\cr
 ~\vdots &&&&& \e &&\bc\ddots\bc&\bc\ddots\bc&&\cr 
  y_{|l|-1} &&&&& \e &&&\bc\ddots\bc&\bc a\bc&\cr
  y_{|l|}    &&&&& \e &&&&\bc b\bc&\bc a\bc\cr}$,
		$Q\ba=\ba\bordermatrix{
	   	&\ba x^+_1\ba&\ba x^+_2\ba&\ba\cdots\ba&\ba x^+_{m-1}\ba&\ba x_m\ba
		&\ba y^+_1\ba&\ba y^+_2\ba&\ba\cdots\ba&\ba y^+_{|l|-1}\bb&\bb y^+_{|l|}\ba\cr
	z_1 	&-b_1\bd&\bd\e_2\bd&&&&&&&&\cr
	z_2 	&&\bd-b_2\bd&\bd\ddots\bd&&&&&&&\cr
   \vdots &&&\bd\ddots\bd&\bd\ddots\bd&&&&&&\cr 
  z_{m-1} &&&&\bd-b_{m-1}\bd&\bd\e_m\bd&&&&&\cr
	z_m &\ba\e_1\ba&&&& \bd-b_m\bd&\e&\e&\e&\e&\e\cr
	w_1 	&&&&&&b&a&&&\cr
	w_2 	&&&&&&&b&a&&\cr
   \vdots &&&&&&&&\bc\ddots\bc&\bc\ddots\bc&\cr 
  w_{|l|-1} &&&&&&&&&\bc b\bc&a\cr
  w_{|l|}    &\ba\e_1\ba&&&&&&&&&b\cr}$	
}	
\pvc
		Then the Alexander polynomial $\De_K(t)$ of $K$ is the product of the Alexander 
		polynomial $\De_k(t)$ of $k$, $|P-t\,Q^T|$, and $|Q-t\,P^T|$. 
%
%------------------------------------------------------------------------------------------------
%
\bgnM
		\label{clm:PQ}
		We have the following, where $c=a-tb$, $d=b-ta$, $e=\e(1-t)$, $c_i=a_i-tb_i$, 
		$d_i=b_i-ta_i$, and $e_i=\e_i(1-t)$.	
		$$|\,Q\ba-\ba tP^T|=\displaystyle 
			d^{\,|l|} \prod_{i=1}^m (-d_i)+(-1)^{|l|+m+1} c^{\,|l|} \prod_{i=1}^m e_i		\,,~\,~\,~
		|P\ba-\ba t\,Q^T|=\displaystyle 
			c^{\,|l|} \prod_{i=1}^m (-c_i)+(-1)^{|l|+m+1} d^{\,|l|} \prod_{i=1}^m e_i.$$
		\endM
%
%------------------------------------------------------------------------------------------------
%
\npg
		{\it Proof.~} First we calculate $|\,P\ba-\ba tQ^T|$ noticing that $e=c+d$. 		
\pva
		If $l=0$, then we have that	$|P\ba-\ba t\,Q^T|=\begin{array}{|ccccc|}
		\ba-c_1\bb&&&&e_1\\[-.5ex]			\ba e_2\ba&\bb-c_2\bc&&&\\[-1ex] 
		&\bc\ddots\bc&\bb\ddots\bc&& \\[-1.5ex] &&\bb\ddots\bc&\ba-c_{m-1}\bc&\\[-.5ex] 
		&&&e_{m}&\bb-c_{m}\bb\end{array}$ 
		$=\displaystyle c^{\,0} \prod_{i=1}^m (-c_i)+(-1)^{0+m+1} d^{\,0} \prod_{i=1}^m e_i$.
\pva
		If $|l|=1$, then we have that
\pva	
		$|P\ba- tQ^T|=\begin{array}{|ccccc|c|}
		\ba-c_1\bb&&&&e_1\bc&e_1\ba\\[-.5ex]\ba 	e_2\ba&\bb-c_2\bc&&&& \\ [-1ex]
		&\bc\ddots\bc&\bb\ddots\bc&&& \\[-1.5ex] 	&&\bb\ddots\bc&\ba-c_{m-1}\bc&& \\[-.5ex] 
		&&&e_{m}&\bb-c_{m}\bb& \\ \hline	&&&&e&c	\end{array}
		=\begin{array}{|ccccc|c|}
		\ba-c_1\bb&&&&0\bc&e_1\ba\\[-.5ex]\ba 	e_2\ba&\bb-c_2\bc&&&& \\[-1ex] 
		&\bc\ddots\bc&\bb\ddots\bc&&& \\[-1.5ex] 	&&\bb\ddots\bc&\ba-c_{m-1}\bc&& \\[-.5ex] 
		&&&e_{m}&\bb-c_{m}\bb& \\ \hline	&&&&d&c	\end{array}\,$
		$=\displaystyle c^{\,1} \prod_{i=1}^m (-c_i)+(-1)^{1+m+1} d^{\,1} \prod_{i=1}^m e_i$.
\pvb
		If $|l|>1$, then we have that
\pva	
		$|P\ba- tQ^T|=\begin{array}{|ccccc|ccccc|}
	\ba-c_1\bb&&&&e_1\bc&&&&&e_1\ba\\[-.5ex]	\ba e_2\ba&\bb-c_2\bc&&&&&&&& \\[-1ex] 
	&\bc\ddots\bc&\bb\ddots\bc&&&&&&&\\[-1.5ex] &&\bb\ddots\bc&\ba-c_{m-1}\bc&&&&&&\\[-.5ex] 
	&&&e_{m}&\bb-c_{m}\bb&&&&& \\ \hline &&&&e&c&&&& \\[-.5ex] &&&&e&d&c&&& \\[-1.2ex] 	
	&&&&e&&\ddots&\ddots&& \\[-1.5ex] &&&&e&&&\ddots&c& \\[-.5ex] &&&&e&&&&d&c\ba\\[-.5ex] 
		\end{array}
		=\begin{array}{|ccccc|ccccc|}
	\ba-c_1\bb&&&&0&&&&&e_1\bb\\[-.5ex]	\ba e_2\ba&\bb-c_2\bc&&&&&&&& \\ [-1ex]
	&\bc\ddots\bc&\bb\ddots\bc&&&&&&&\\[-1.5ex] &&\bb\ddots\bc&\ba-c_{m-1}\bc&&&&&&\\[-.5ex] 
	&&&e_{m}&\bb-c_{m}\bb&&&&& \\ \hline &&&&d&c&&&& \\[-.5ex] &&&&0&d&c&&& \\[-1.2ex] 	
	&&&&0&&\ddots&\ddots&& \\[-1.5ex] &&&&0&&&\ddots&c&\\[-.5ex] &&&&0&&&&d&c\bb\\[-.5ex] 
		\end{array}\,$
\pvan\hsp{20mm}
		$=\displaystyle c^{\,|l|} \prod_{i=1}^m (-c_i)+(-1)^{|l|+m+1} d^{\,|l|} \prod_{i=1}^m e_i$.
		
%
%------------------------------------------------------------------------------------------------
%
\pvc		
		Next we calculate $|\,Q\ba-\ba tP^T|$ noticing that $e=c+d$. 
\pva
		If $l=0$, then we have that 
		$|Q-t\,P^T|=\begin{array}{|ccccc|} \bb-d_1\bb&e_2\bc&&& \\[-1ex]
		&\bb-d_2\bb\bc&\ddots\bb&& \\[-1.2ex] 	&&\ddots\bb&\ddots\bb&\\[-.5ex] 
		&&&\bd-d_{m-1}\bc&e_{m} \\[-.5ex] 		e_1\bc&&&&\bc-d_{m}\bb\end{array}
		=\displaystyle d^{\,0} \prod_{i=1}^m (-d_i)+(-1)^{0+m+1} c^{\,0} \prod_{i=1}^m e_i$.
\pva
		If $|l|=1$, then we have that
\pva
		$|\,Q\ba- tP^T|=\begin{array}{|ccccc|c|}
		\bb-d_1\bb&e_2\bc&&&& \\[-.5ex]	&\bb-d_2\bb\bc&\ddots\bb&&& \\[-1.2ex] 
		&&\ddots\bb&\ddots\bb&& \\[-.5ex] &&&\bd-d_{m-1}\bc&e_{m}& \\[-.3ex] 
		e_1\bc&&&&\bc-d_{m}\bb&e \\ \hline	e_1\bc&&&&&d \\[-.5ex]	\end{array}
		=\begin{array}{|ccccc|c|}
		\bb-d_1\bb&e_2\bc&&&& \\[-.5ex]	 &\bb-d_2\bb\bc&\ddots\bb&&& \\[-1.2ex]
		&&\ddots\bb&\ddots\bb&& \\[-.5ex] &&&\bd-d_{m-1}\bc&e_{m}& \\[-.3ex] 
		0&&&&\bc-d_{m}\bb&c\\ \hline	e_1\bc&&&&&d \\[-.5ex]		\end{array}$
		$=\displaystyle d^{\,1} \prod_{i=1}^m (-d_i)+(-1)^{1+m+1} c^{\,1} \prod_{i=1}^m e_i$.
\pvb	
		If $|l|>1$, then we have that
\pva
		$|\,Q\ba- tP^T|=\begin{array}{|ccccc|ccccc|}
	\bb-d_1\bb&e_2\bc&&&&&&&& \\[-.5ex]	&\bb-d_2\bb\bc&\ddots\bb&&&&&&& \\ [-1.2ex]
	&&\ddots\bb&\ddots\bb&&&&&& \\[-.5ex] 	&&&\bb-d_{m-1}\bc&e_{m}&&&&& \\[-.3ex] 
	e_1\bc&&&&\bc-d_{m}\bb&e&e&\cdots&e&e \\ \hline
	&&&&&d&c&&& \\[-.5ex] &&&&&&d&\ddots&& \\[-1.2ex] &&&&&&&\ddots&\ddots& \\[-.5ex]
	&&&&&&&&d&c \\[-.5ex] e_1\bc&&&&&&&&&d \\[-.5ex] 
		\end{array}
		=\begin{array}{|ccccc|ccccc|}
	\bb-d_1\bb&e_2\bc&&&&&&&& \\[-.5ex]	&\bb-d_2\bb\bc&\ddots\bb&&&&&&& \\[-1.2ex]
	&&\ddots\bb&\ddots\bb&&&&&& \\[-.5ex] 	&&&\bb-d_{m-1}\bc&e_{m}&&&&& \\[-.3ex] 
	0&&&&\bc-d_{m}\bb&c&0&\cdots&0&0\\ \hline
	&&&&&d&c&&& \\[-.5ex] &&&&&&d&\ddots&& \\[-1.2ex] &&&&&&&\ddots&\ddots& \\[-.5ex]
	&&&&&&&&d&c \\[-.5ex] e_1\bc&&&&&&&&&d \\[-.5ex] 
		\end{array}$
\pvan\hsp{20mm}
		$=\displaystyle d^{\,|l|} \prod_{i=1}^m (-d_i)+(-1)^{|l|+m+1} c^{\,|l|} \prod_{i=1}^m e_i$.
\vsme	\EndR
%
%------------------------------------------------------------------------------------------------
%
\npg
		Now we calculate the Alexander polynomial $\De_K(t)$ of $K$ diving the case into two 
		depending on the value of $l$; $l\geq 0$ or $l<0$. Here note the following.
		$${\renewcommand\arraystretch{1.2}
		\begin{array}{|c||c|c|c|c|c|} \hline\e&a&b&c&d&e\\ 
		\hline 1&1&0&1&-t&1-t\\ 	\hline -1&0&-1&t&-1&-(1-t)\\\hline\end{array}~\;~\;~\;~
		\begin{array}{|c||c|c|c|c|c|} \hline\e_i&a_i&b_i&c_i&d_i&e_i\\ 
		\hline 1&1&0&1&-t&1-t\\ 	\hline -1&0&-1&t&-1&-(1-t)\\\hline\end{array}}$$
\pvcn
		{\bf Case} $l\geq0$ : From the above table, we have the following;\pvan
		$|P-t\,Q^T|=\displaystyle c^{\,|l|} \prod_{i=1}^m (-c_i)+(-1)^{|l|+m+1} d^{\,|l|} \prod_{i=1}^m e_i$
\pcg\hse $=1^{\, l} \,(-1)^p (-t)^{m-p} + (-1)^{l+m+1} (-t)^{\, l} (-1)^{m-p} (1-t)^m$
		$=(-1)^{1-p}\,\{t^{\, l} (1-t)^m-(-t)^{m-p}\}$		
\pvcn		
		$|Q-t\,P^T|=\displaystyle d^{\,|l|} \prod_{i=1}^m (-d_i)+(-1)^{|l|+m+1} c^{\,|l|} \prod_{i=1}^m e_i$
\pcg\hse $=(-t)^{\,l}\,t^{\,p}\, 1^{m-p} + (-1)^{l+m+1}1^{\, l} \, (-1)^{m-p} (1-t)^m$
		$=(-1)^{l+1-p}\,\{(1-t)^m-t^{\,l}(-t)^p\}$
\pvcn
		{\bf Case} $l<0$ : From the above table, we have the following;\pvan
		$|P-t\,Q^T|=\displaystyle c^{\,|l|} \prod_{i=1}^m (-c_i)+(-1)^{|l|+m+1} d^{\,|l|} \prod_{i=1}^m e_i$
\pcg\hse $=t^{-l} \,(-1)^p (-t)^{m-p} + (-1)^{-l+m+1} (-1)^{-l} (-1)^{m-p} (1-t)^m$
		$=(-1)^{1-p}\,\{(1-t)^m-t^{-l}(-t)^{m-p}\}$
\pvcn		
		$|Q-t\,P^T|=\displaystyle d^{\,|l|} \prod_{i=1}^m (-d_i)+(-1)^{|l|+m+1} c^{\,|l|} \prod_{i=1}^m e_i$
\pcg\hse $=(-1)^{-l} \,t^{\,p}\, 1^{m-p} + (-1)^{-l+m+1}t^{-l} \, (-1)^{m-p} (1-t)^m$
		$=(-1)^{-l+1-p}\,\{t^{-l}(1-t)^m-(-t)^p\}$				
\pvc
		In both cases, we obtain that 
		$\De_K(t)\doteq\De_k(t)\,\{(1-t)^m - t^{\,l} (-t)^p\}\,\{(1-t)^m - t^{-l} (-t)^{m-p}\}$,
		and thus we complete the proof.
		\EndR
\pvc
\bgnF  		
		{\iclg[scale=.9, bb=0 0 357 176]{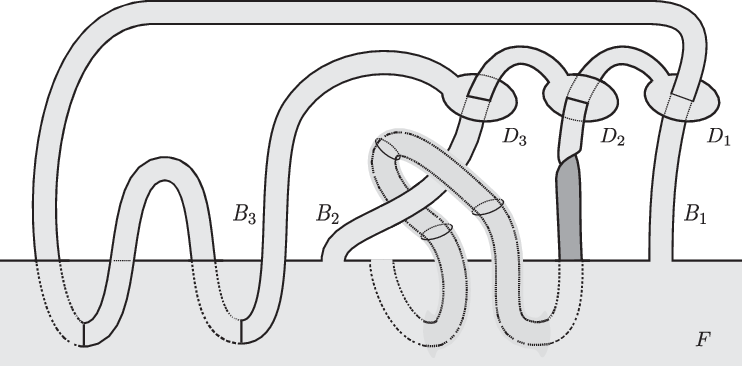}} \captn{}\label{fig:stp_expl}	\endF
\npg
\bgnF  		
		{\iclg[scale=.9, bb=0 0 358 194]{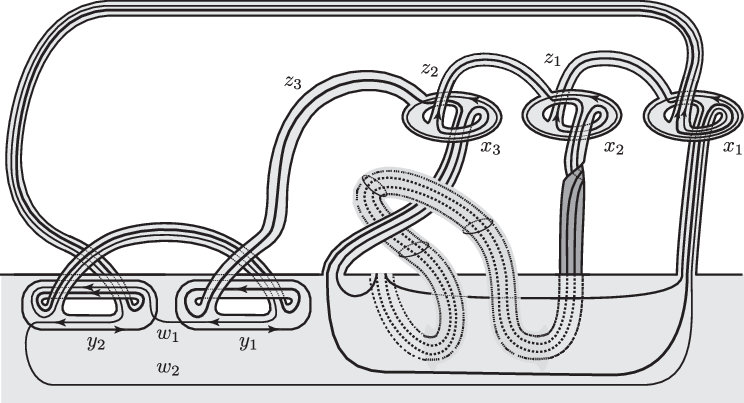}} \captn{}\label{fig:exp_basis}\endF			
	
%-------------------------------------------------------------------------------------------------------------------------------------------------	
%
%
%
%
%
%
%
%
%
%-------------------------------------------------------------------------------------------------------------------------------------------------	
\pvcn
		{\bf\it Proof of Theorem \ref{thm:nsrk}.}	
		For each $i$ $(1\leq i \leq N)$, we can construct a simple-ribbon knot $k_i$ with
		$\De_{k_i}(t)=\v(t\,;m_i,p_i,l_i)\, \v(t^{-1}\,;m_i,p_i,l_i)$ by following the proof of Theorem
		\ref{thm:alexpoly} (see also Figure \ref{fig:stp_expl}). Let $K^\ast$ be the connected 
		sum of $k_1$, $k_2$, $\ldots$, $k_N$. Then $K^\ast$ is a simple-ribbon knot such 
		that $\De_{K^\ast}(t)=\De(t)$. Let $\mD\cup\mB$ be the set of disks and bands 
		which gives the SR-fusion on the trivial knot $\mO=\ptl D_0$ producing $K^\ast$. 
		Take a $3$-ball $H$ which is a small neighborhood of a point of $\mO-\mB$ and
		a trivial knot $\r$ in $H$ which intersects $D_0$ twice so that $\lk(\r, \mO)=2$.
		Let $V^\ast$ be the closure of $S^3-N(\r; S^3)$. Since $\r$ is the trivial knot,
		$V^\ast$ is an unknotted torus which contains $K^\ast$ with $w_{V^\ast}(K^\ast)=2$,
		where $w_{V^\ast}(K^\ast)$ is the absolute value of the algebraic intersection 
		number of $K^\ast$ with a meridian disk of $V^\ast$.
\pvc
		Let $V$ be a tubular neighborhood of the Kinoshita-Terasaka knot $\k$ and 
		$f$ a faithful homeomorphism of $V^\ast$ onto $V$, i.e. $f$ maps the preferred 
		longitude of $\ptl V^\ast$ onto the preferred longitude of $\ptl V$. 
		Since $\De_\k(t)=1$, we obtain that $\De_K(t)=\De_{K^\ast}(t)\,\De_\k(t^2)$
		$=\De_{K^\ast}(t)=\De(t)$ for $K=f(K^\ast)$ by Proposition 8.23 of \cite{BZ-85}.
		Since $f$ is faithful and both of $K^\ast$ and $\k$ are ribbon knots,
		$K$ is also a ribbon knot by Lemma 3 of \cite{TS-MSNK80}\footnote{
		Lemma 3 shows that $K$ is ribbon cobordant to $K^\ast$ if $\k$ is a ribbon knot,
		although it states that $K$ is cobordant to $K^\ast$.}.
		On the other hand, as $w_{V}(K)=w_{V^\ast}(K^\ast)=2$ and $\k$ is a non-trivial
		knot, $K$ is not a simple ribbon knot by Corollary 1.8 of \cite{KST-SRCP18}.
		\EndR

%-------------------------------------------------------------------------------------------------------------------------------------------------
%
%
%
%
%
%
%
%
%
%
%
%
%
%
%
%
%
%
%
%
%
%
%
%
%
%
%
%
%
%
%-------------------------------------------------------------------------------------------------------------------------------------------------

					\section{Proof of Theorem \ref{thm:msrk}}\label{sec:thmmsrk}
					
		Note that if $K$ is a knot of $\mK_m$, then $\det(K)=|\De_K(-1)|=(2^{\,m}-1)^a (2^{\,m}+1)^b$
		for some non-negative integers $a$ and $b$ by Corollary \ref{cor:alexpoly}. 
		Moreover if $K$ is also a knot of $\mK_n$, then $\det(K)=(2^{\,n}-1)^c (2^{\,n}+1)^d$ 
		for some non-negative integers $c$ and $d$, and thus the set of prime factors of
		$(2^{\,m}-1)^{a'}(2^{\,m}+1)^{b'}$ and $(2^{\,n}-1)^{c'}(2^{\,n}+1)^{d'}$ coinside,
		where $i'={\rm min}(i, 1)$ for $a$, $b$, $c$, and $d$.
\pva
		Let $P(x)$ be the set of prime factors of an integer $x>1$, and $(y,z)$ the greatest 
		common divisor of positive integers $y$ and $z$. Note that if $P(y)=P(z)$ and 
		$(y, z)=w$, then we have that $P(y)=P(z)=P(w)$. Here we prepare several lemmas, 
		the first one of which is the theorem by P. Mih\u{a}ilescu (the Catalan conjecture).
		
\bgnL	{\rm(\cite[Theorem~5]{PM-JRAM04})}\label{lem:cc}	
		The following equation has no other integer solutions but $3^2-2^3=1$.
%		There are no solutions of the Diophantine equation
		\begin{equation}\label{eqcc} x^u-y^v=1~(x>0, y>0, u>1, v>1) \end{equation}
%		other than $x^u=3^2$, $y^v=2^3$.
		\endL
		
\bgnL 	{\rm(\cite[Theorem~1]{TI-MAAM04})}\label{lem:mm}
		Let $A$, $m$, and $n$ be integers such that $A>1$ and $m>n\geq 1$.
		Then $P(A^m-1)=P(A^n-1)$ if and only if $m=2$, $n=1$,  and $A=2^l-1$
		for an integer $l>0$.
		\endL

\bgnL	\label{lem:ppm} Let $A$ be an integer such that $A>1$. Then the followings hold.
		\bgnI	
	\item[(1)] $P(A^p+1)=P(A+1)$ for an odd integer $p$ $(>1)$ if and only if $p=3$ and $A=2$.
	\item[(2)] $P(A^q -1)=P(A+1)$ for an even integer $q$ $(>0)$ if and only if $q=2$ and 
			$A=2^l+1$ for an integer $l\geq 0$.
		\endI\endLR
		Since the if parts are easy to be checked, we only show the only if parts.
\pvan
	(1) 	First the following equation holds, since $p$ is odd.
		\begin{equation}\label{eq3a}
		B=\Frc{A^p+1}{A+1}
		  =A^{p-1}-A^{p-2}+\cdots-A+1=\Sum_{i=0}^{p-2} \binom{p}{i} (A+1)^{p-i-1} (-1)^i +p		
		\end{equation}

		If $p$ is prime, then we have that $(B,A+1)=(A+1,p)=p$ from equation (\ref{eq3a}),
		and thus that $P(B)=\{p\}$, since $P(B)\subset P(A^p+1)=P(A+1)$. Moreover, we have 
		that $B\equiv p$ $($mod $p^2)$ also from equation (\ref{eq3a}), since $A+1\equiv 0$ 
		$($mod $p)$, $\binom{p}{p-2}\equiv 0$ $($mod $p)$. Hence we obtain that $B=p$. 
		If $p>3$, then we also have that \begin{equation}\label{eq3b}
		B=A^{p-1}-A^{p-2}+\cdots-A+1=A(A-1)(A^{p-3}+A^{p-5}+\cdots+1)+1>A(A-1)\Frc{p-1}{2}+1
		\geq p,
		\end{equation}
		since $A\geq 2$. However then it contradicts that $B=p$. Therefore we have 
		that $p=3$. Then we have that $A^2-A+1=B=p=3$ from equation (\ref{eq3a}), 
		and thus that $A=2$, since $A>1$, which completes the proof.
\pvb
		If $p$ is not prime, then let $p'$ be a prime factor of $p$, and let $p=p'r$ and
		$B=A^r$. Since $r$ and $p'$ are odd, we have that $A+1$ divides $A^r+1=B+1$
		and that 	$B+1$ divides $B^{p'}+1$. Hence we have that 
		$P(A+1)\subset P(B+1)\subset P(B^{p'}+1)=P(A^p+1)$, since $B^{p'}=A^p$.
		Hence we have that $P(B^{p'}+1)=P(B+1)$, since $P(A^p+1)=P(A+1)$. 
		Thus from the previous case, we have that $p'=3$ and $B=A^r=2$, and thus 
		$A=2$ and $r=1$. However then, we have that $p=p'r=3$, which contradicts
		that $p$ is not prime.
\pvbn
	(2)	Since $q$ is even, we have that $q\geq 2$. Hence we have that $P(A-1)\subset 
		P(A^q-1)=P(A+1)$, and thus that $P(A^2-1)=P((A-1)(A+1))=P(A+1)=P(A^q-1)$. 
		Thus we have that $q=2$ from Lemma \ref{lem:mm}. If $A\neq 2=2^0+1$, then 
		we have that $A-1>1$ and thus that $A+1$ and $A-1$ 
		are not coprime, since $P(A-1)\subset P(A+1)$. Hence we have that 
		$(A+1, A-1)=(A-1, 2)=2$, since $A+1=(A-1)+2$. Therefore we obtain that 
		$A-1=2^l$ for $l>0$, which completes the proof.
		\endR			
		
		Using Lemma \ref{lem:cc} and Lemma \ref{lem:ppm}, we show the following.
	
\bgnP	\label{prop:key}
		Let $A$, $m$, and $n$ be integers such that $A>1$ and $m$, $n\geq 1$.
		Then we have the following.\bgnI
		\item[(1)] 	$P(A^m+1)=P(A^n+1)$ $(m>n)$ if and only if $m=3$, $n=1$, and $A=2$;\pva
		\item[(2)] 	$P(A^m+1)=P(A^n-1)$ if and only if one of the following holds.\bgnI
		\item[(i)]		$m=1$, $n=1$, and $A=3\,;$
		\item[(ii)]		$m=3$, $n=2$, and $A=2\,;$
		\item[(iii)]	$m=2$, $n=4$, and $A=3\,;$ and
		\item[(iv)]	$m=1$, $n=2$, and $A=2^l+1$ 	for an integer $l\geq 0$.
		\endI\endI\endPR
		First we have the following for indeterminate $X$ and positive integers 
		$s$, $t$, and $q$	and a non-negative integer $r$ such that $s=qt+r$.
	\begin{equation}\label{eqp1}
		X^s\pm 1=(X^t+1)(X^{s-t}-X^{s-2t}+\cdots-(-1)^qX^{s-qt})+(-1)^qX^r\pm 1
		\end{equation}
	\begin{equation}\label{eqp2}
		X^s+1=(X^t-1)(X^{s-t}+X^{s-2t}+\cdots+X^{s-qt})+X^r+1~\ \ \ \ \ \ \ ~
		\end{equation}

		Let $g=(m,n)$. Then we have the following.

\bgnM	\label{clm:pp}	 $(A^m+1, A^n+1), (A^m+1, A^n-1)=1$, $2$ or $A^g+1$.
\endMR
		For positive integers $c_0$ and $c_1$, 
		let $(c_0, c_1)=(c_1,c_2)=\cdots=(c_{N-1},c_N)=c_N$ be the sequence obtained by 
		the Euclidian algorithm. Then letting $c_i=q_{i+1}c_{i+1}+q_{i+2}$, we also have the 
		following from equations (\ref{eqp1}) and (\ref{eqp2}).
	\begin{equation}\label{eqcpp1}
		A^{c_{N-1}}\pm 1=(A^{c_N}+1)(A^{c_{N-1}-c_N}-A^{c_{N-1}-2c_N}+\cdots
								-(-1)^qA^{c_{N-1}-{q_N}c_N})+(-1)^{q_N}A^0\pm 1
		\end{equation}
	\begin{equation}\label{eqcpp2}
		A^{c_{N-1}}+1=(A^{c_N}-1)(A^{c_{N-1}-c_N}+A^{c_{N-1}-2c_N}+\cdots
								+A^{c_{N-1}-{q_N}c_N})	+A^0+1~\ \ \ \ \ \ \ \ ~ \
		\end{equation}
		Hence by letting $(c_0, c_1)=(m,n)$ or $(n,m)$, we have that
		$(A^m+1, A^n+1)$, $(A^m+1, A^n-1)$ is either $A^g+1$ or $(A^g\pm 1, 2)$, 
		which induces the conclusion.
\endR
\pva
		Since the if parts are easy to be checked, we only show the only if parts.
\pvcn
	(1)	Since $P(A^m+1)=P(A^n+1)$, we have that $A^m+1$ and $A^n+1$ are not coprime, 
		and thus that $(A^m+1, A^n+1)=2$ or $A^g+1$ from Claim \ref{clm:pp}.
		In the former case, we have that $P(A^m+1)=P(A^n+1)=P(2)=\{2\}$.
		Thus, $A^m+1=2^k$ and $A^n+1=2$ for $k>1$, since $m>n$.
		However then, we have that $A=1$, which contradicts that $A>1$.
		In the latter case, we have that $P(A^m+1)=P(A^n+1)=P(A^g+1)$ and that
		$m=gM$ with an odd integer $M$ from equation (\ref{eqp1}).
		If $M=1$, then $m=g$, which contradicts that $m>n$. Thus $M$ is odd and $M>1$.  
		Then we have that $M=3$ and $A^g=2$ by Lemma \ref{lem:ppm} (1), and thus that
		$A=2$, $g=1$, $m=gM=3$. Hence we have that $n=g=1$, since $m>n$,
		which completes the proof. 
\pvcn	
	(2)	Since $P(A^m+1)=P(A^n-1)$, we have that $A^m+1$ and $A^n-1$ are not coprime, 
		and thus that $(A^m+1, A^n-1)=2$ or $A^g+1$ from Claim \ref{clm:pp}.
		In the former case, we have that $P(A^m+1)=P(A^n-1)=P(2)=\{2\}$, and thus 
		that $A^m+1=2$ or $A^n-1=2$. If $A^m+1=2$, then $A^m=1$, which contradicts 
		that $A>1$. If $A^m+1=2^k$ $(k>1)$ and $A^n-1=2$, then we have that $A=3$ 
		and $n=1$, and thus that $A^m+1=3^m+1=2^k$ $(k>1)$. Then by Lemma \ref{lem:cc},
		we have that $m=1$, and thus obtain condition (i).
\pvc		
		In the latter case, we have that $P(A^m+1)=P(A^n-1)=P(A^g+1)$ and that $m=gM$ 
		with an odd integer $M$ from equation (\ref{eqp1}). Consider the case where $M>1$. 
		Then we have that $M=3$ and $A^g=2$ by Lemma \ref{lem:ppm} (1), and thus that 
		$A=2$, $g=1$, $m=gM=3$. Since $A^m+1=2^3+1=9$, and thus $P(2^n-1)=P(9)=\{3\}$ 
		and $(A^m+1, A^n-1)=(9,2^n-1)=3$, we have that $2^n-1=3$ and thus that $n=2$. 
		Therefore we obtain condition (ii).
\pva	
		Next consider the case where $M=1$, i.e., $m=g$. Hence let $n=mq$ $(q\geq 1)$
		and $D=A^m$. Thus we have that $P(D+1)=P(D^q-1)$ and that $(D+1, D^q-1)=D+1$.
		Therefore $q$ is even, since otherwise $D+1$ does not divide $D^q-1$. Then we 
		have that $q=2$ and $D=2^l+1$ for $l\geq 0$ by Lemma \ref{lem:ppm} (2). 
		If $m>1$ and $l>1$, then the equation $A^m=2^l+1$ has the unique solution 
		$(A, m, l)=(3,2,3)$ by Lemma \ref{lem:cc}, and thus we obtain condition (iii). 
		If $m=1$, then we have that $n=mq=2$ and $A=2^l+1$ for $l\geq 0$, i.e., 
		condition condition (iv). If $l=0$ (resp. $1$), then we have that $A=2$ (resp. $A=3$)
		and $m=1$, and thus that condition (iv).
		\endR
		
		Now using Proposition \ref{prop:key} and Lemma \ref{lem:mm}
		we obtain the following.
		
\bgnL	\label{lem:main}
		Let $p$, $q$, $r$, $s$, $M$, $N$ be positive integers with $M\neq N$.
		Then we have the following.
\bgnI	\item[(1)]  	$(2^M-1)^p \neq (2^N-1)^r$. \pva
		\item[(2)]  If	$(2^M+1)^q=(2^N+1)^s$ $(M>N)$, then $M=3$, $N=1$, and $s=2q$.\pva
		\item[(3)]  If 	$(2^M+1)^q=(2^N -1)^r$, then 
					$M=3$, $N=2$, $r=2q$ or $M=1$, $N=2$, $q=r$.\pva
		\item[(4)]	$(2^M-1)^p (2^M+1)^q \neq (2^N-1)^r (2^N+1)^s$\pva
		\item[(5)]  If  $(2^M-1)^p (2^M+1)^q = (2^N-1)^r$, then
					$2M=N$, $p=q=r$.\pva
		\item[(6)]  If	$(2^M-1)^p (2^M+1)^q = (2^N+1)^r$, then
					$M=1$, $N=3$, $q=2r$.
		\endI\endLR
		Note that if positive integers $X$, $Y$ and non-negative integers 
		$p$, $q$ satisfies the equation $X^p=Y^q$, then $P(X)=P(Y)$. 
		The first three statements are obtained by Lemma \ref{lem:mm}, 
		Proposition \ref{prop:key} (1), and Proposition \ref{prop:key} (2), respectively.
		For the last three statements, note that $P((2^M-1)^p (2^M+1)^q)=P(2^{2M}-1)$. 
		Therefore (4) and (5) are obtained by Lemma \ref{lem:mm}, and
		(6) is obtained by Proposition \ref{prop:key} (2).
		\endR
\pvcn
		{\bf\it Proof of Theorem \ref{thm:msrk}.}	
		Let $K$ be a knot of $\mK_m\cap\mK_n$. Then we have that 
		$\det(K)=(2^{\,m}-1)^a (2^{\,m}+1)^b=(2^{\,n}-1)^c (2^{\,n}+1)^d$ for some 
		non-negative integers $a$, $b$, $c$, and $d$ by Corollary \ref{cor:alexpoly}.
		Thus we obtain the conclusion by Lemma \ref{lem:main}. 
		\EndR

{\small
\renewcommand{\arraystretch}{1.5}
\begin{table}[htb]
\caption{Ribbon knots with up to 10 crossings, where $F(t\,;m,p,l)=\v(t\,;m,p,l)\,\v(t^{-1}\,;m,p,l)$}\label{tab:alex}
  \begin{tabular}{|c||c|c|c|r l|} \hline
    $K$ & simple-ribbon & $\d_2(K)$ & det(K) & & $\Delta'_K(t)$   \\ \hline\hline
    $6_1$ & $\bigcirc$ & 0 & 9 		& $F(t;2,0,0)=$&$2\ms 5t+ 2t^2$  \\ \hline
    $8_8$ & $\bigcirc$ & 5 & 25 	& $F(t;2,1,-1)=$&$2\ms 6t+ 9t^2\ms 6t^3+ 2t^4$     \\ \hline
    $8_9$ & $\bigcirc$ & 7 & 25 	& $F(t;2,2,1)=$&$1\ms 3t+ 5t^2\ms 7t^3+ 5t^4\ms 3t^5+ t^6$    \\ \hline
    $8_{20}$ & $\bigcirc$ & 9 & 9 	& $F(t;2,1,0)=$&$1\ms 2t+ 3t^2\ms 2t^3+ t^4$ \\ \hline
    $9_{27}$ & $\bigcirc$ & 5 & 49 	& $F(t;3,1,0)=$&$1\ms 5t+ 11t^2\ms 15t^3+ 11t^4\ms 5t^5+ t^6$\\ \hline
    $9_{41}$ & $\bigcirc$ & 7 & 49	& $F(t;3,0,0)=$&$3\ms 12t+ 19t^2\ms 12t^3+ 3t^4$  \\ \hline
    $9_{46}$ & $\bigcirc$ & 0 & 9	& $F(t;2,0,0)=$&$2\ms 5t+ 2t^2$ 	\\ \hline
    $10_{3}$ & $\times$ & 1 & 25	& &$6\ms 13t+ 6t^2$   \\ \hline
    $10_{22}$ & $\times$ & 11 & 49	& &$2\ms 6t+ 10t^2\ms 13t^3+ 10t^4\ms 6t^5+ 2t^6$      \\ \hline
    $10_{35}$ & $\times$ & 1 & 49 	& &$2\ms 12t+ 21t^2\ms 12t^3+ 2t^4$     \\ \hline
    $10_{42}$ & $\bigcirc$ & 9 & 81	& $F(t;3,2,-1)=$&$1\ms 7t+ 19t^2\ms 27t^3+ 19t^4\ms 7t^5+ t^6$     \\ \hline
    $10_{48}$ & $\times$ & 91 & 49	& &$1\ms 3t+ 6t^2\ms 9t^3+ 11t^4\ms 9t^5+ 6t^6\ms 3t^7+ t^8$    \\ \hline
    $10_{75}$ & $\bigcirc$ & 9 & 81	& $F(t;3,3,-1)=$&$1\ms 7t+ 19t^2\ms 27t^3+ 19t^4\ms 7t^5+ t^6$      \\ \hline
    $10_{87}$ & $\bigcirc$ & 0 & 81	& $F(t;3,2,-1)=$&$2\ms 9t+ 18t^2\ms 23t^3+ 18t^4\ms 9t^5+ 2t^6$    \\ \hline
    $10_{99}$ & $\bigcirc$ & 81 & 81 & 
    				$F(t;1,1,1) F(t;2,1,0)=$&$1\ms 4t+ 10t^2\ms 16t^3+ 19t^4\ms 16t^5+ 10t^6\ms 4t^7+ t^8$ \\ \hline
    $10_{123}$ & $\times$ & 1 & 121 	& &$1\ms 6t+ 15t^2\ms 24t^3+ 29t^4\ms 24t^5+ 15t^6\ms 6t^7+ t^8$    \\ \hline
    $10_{129}$ & $\bigcirc$ & 5 & 25 	& $F(t;2,1,1)=$&$2\ms 6t+ 9t^2\ms 6t^3+ 2t^4$    \\ \hline
    $10_{137}$ & $\bigcirc$ & 1 & 25	& $F(t;2,0,1)=$&$1\ms 6t+ 11t^2\ms 6t^3+ t^4$      \\ \hline
    $10_{140}$ & $\bigcirc$ & 9 & 9		& $F(t;2,1,0)=$&$1\ms 2t+ 3t^2\ms 2t^3+ t^4$      \\ \hline
    $10_{153}$ & $\bigcirc$ & 35 & 1	& $F(t;1,1,2)=$&$1\ms t\ms t^2+ 3t^3\ms t^4\ms t^5+ t^6$     \\ \hline
    $10_{155}$ & $\bigcirc$ & 7 & 25	& $F(t;2,2,1)=$&$1\ms 3t+ 5t^2\ms 7t^3+ 5t^4\ms 3t^5+ t^6$ \\ \hline
    $3_1\sh 3_1^*$ & $\bigcirc$ & 9 & 9 	& $F(t;1,1,1)=$&$1\ms 2 t+3 t^2\ms 2t^3+t^4$     \\ \hline
    $4_1\sh 4_1$ & $\bigcirc$ & 1 & 25 	& $F(t;2,0,1)=$&$1\ms 6 t+11 t^2\ms 6 t^3+t^4$     \\ \hline
    $5_1\sh 5_1^*$ & $\times$ & 121 & 25  	& &$1\ms 2 t+3 t^2\ms 4 t^3+5 t^4\ms 4 t^5+3 t^6\ms 2 t^7+t^8$   \\ \hline
    $5_2\sh 5_2^*$ & $\times$ & 1 & 49& &$4 \ms 12 t+17 t^2\ms 12 t^3+4t^4$  \\ \hline
  \end{tabular}
\end{table}
}	
%-------------------------------------------------------------------------------------------------------------------------------------------------
						\section*{Acknowledgement} 
          	
		The authors would like to thank Alexander Zupan for informing us errors of the diagrams of
		$10_{42}$ and $10_{75}$ in Figure \ref{fig:exp_srkten} of \cite{KSTI-OJM21}.
		
\vsp{-1mm}
 						
{\small
\pvcn	
		Kengo KISHIMOTO (e-mail:kengo.kishimoto@oit.ac.jp) \\
		Tetsuo SHIBUYA    \\
		Tatsuya TSUKAMOTO (e-mail:tatsuya.tsukamoto@oit.ac.jp)\\
		Tsuneo ISHIKAWA  (e-mail:tsuneo.ishikawa@oit.ac.jp) \\
		Department of Mathematics,		Osaka Institute of Technology,	
		Asahi,Osaka 535-8585, Japan \\
}

\begin{thebibliography}{99}
						
\bibitem{BZ-85}            		G.~Burde and H.~Zieschang,
                       			   	Knots, 
                               			de Gruyter, 1985.
	
\bibitem{TI-MAAM04}	T.~Ishikawa, N.~Ishida and Y.~Yukimoto,
							{\em On prime factors of $A^n-1$},
							Amer. Math. Monthly, {\bf 111} (2004), 243--245.
	
\bibitem{AK-book96}	 	A.~Kawauchi,
							A survey of knot theory,
							Birkhuser Verlag, Basel, 1996.
											
\bibitem{KST-JMSJ16}	 K.~Kishimoto, T.~Shibuya and T.~Tsukamoto,
							{\em Simple-ribbon fusions and genera of links},
							J. Math. Soc. Japan, {\bf 68} (2016), 1033--1045.	
							
\bibitem{KST-SRCP18}	 K.~Kishimoto, T.~Shibuya and T.~Tsukamoto,
							{\em Simple-ribbon concordance of knots},
							Kobe J. Math, {\bf 37} (2020), 1--17.
							
\bibitem{KSTI-OJM21}	 K.~Kishimoto, T.~Shibuya T.~Tsukamoto and T.~Ishikawa,
							{\em Alexander polynomials of simple-ribbon knots},
							Osaka J. Math, {\bf 58} (2021), 41--57.							
						
\bibitem{PM-JRAM04}	P.~Mih\u{a}ilescu,
							{\em Primary Cyclotomic Units and a Proof of Catalan's Conjecture},
							J. Reine Angew. Math., {\bf 572} (2004), 167--195.
						
\bibitem{HS-ACTM53} 	H.~Schubert,
							{\em Knoten und vollringe},
							Acta Mathematica, {\bf 90} (1953), 131--286.	
						
\bibitem{TS-MSNK80}	T.~Shibuya,
							{\em On the cobordism of compound knots},
							Math. Sem. Notes, Kobe Univ., {\bf 8} (1980), 331--337.						

          		\end{thebibliography}
\end{document}